\documentclass[a4paper,12pt,leqno]{article}
\usepackage{graphicx}
\usepackage{amsmath,amssymb,amsthm}
\usepackage{mathrsfs}
\numberwithin{equation}{section}
\newtheorem{theorem}{Theorem}[section]
\newtheorem{proposition}{Proposition}[section]

\title{\textbf{Hypersurfaces of a Projective Randers conformal change}}
\date{}
\begin{document}
\maketitle
\begin{center}
\author{\textbf{V. K. Chaubey and $^{*}$Pradeep Kumar}}
\end{center}
\begin{center}
Department of Applied Sciences, Buddha Institute of Technology \\* Sector-7, Gida, Gorakhpur (U.P.)-273209, INDIA, E-Mail: vkchaubey@outlook.com
\end{center}
\begin{center}
$ ^{*} $Department of Mathematics and Statistics, DDU Gorakhpur University \\ Gorakhpur (U.P.)-273209, INDIA
\end{center}
\begin{abstract} 
\hspace{10pt} In the year 1984 Shibata investigated the theory of a change which is called a $ \beta $-change of a Finsler metric. On the other hand in 1985 a systematic study of geometry of hypersurfaces in Finsler spaces was given by Matsumoto. In the present paper is to devoted to the study of a condition for a Randers conformal change to be projective and find out when a totally geodesic hypersurface $ F^{n-1} $ remains to be a totally geodesic hypersurface $ F^{n-1} $ under the projective Randers conformal change. Further obatined the condition under which a Finslerian hypersurfaces given by the projective Randers conformal change are projectively flat. \newline \\*
\end{abstract}
\bigskip
\section{Introduction}
\hspace{10pt} Let $ (M^{n}, L) $ be an n-dimensional Finsler space on a differentiable manifold $ M^{n} $, equipped with the fundamental function $ L(x, y) $. In 1984, Shibata \cite{Sh} introduced the transformation of Finsler metric:
\begin{equation}
L^{'}(x, y)=f(L, \beta)
\end{equation}
where $ \beta=b_{i}(x)y^{i}, \;\; b_{i}(x) $ are components of a covariant vector in $ (M^{n}, L) $ and $ f $ is positively homogeneous function of degree one in $ L $ and $ \beta $. This change of metric is called a $ \beta $-change. \newline \\*
\hspace{10pt} The conformal theory of Finsler spaces has been initiated by M.S. Knebelman \cite{Kne} in 1929 and has been investigated in detail by many authors \cite{Ha, Iz1, Iz2, Kne} etc. The conformal change is defined as
\begin{center}
$ L(x,y)\rightarrow e^{\sigma(x)}L(x, y) $,
\end{center}
where $ \sigma(x) $ is a function of position only and known as conformal factor.\newline \\*
On the other hand in 1985 M. Matsumoto investigated the theory of Finslerian hypersurface \cite{M2}. He has defined three types of hypersurfaces that were called a hyperplane of the first, second and third kinds. \newline\\*
Finslerian Hypersurfaces for a change in a Fisnler metric was studied by several authors \cite{K1, K2, K3, PT, SCM}and obtained so many important results in the stand point of Finsler geometry.\newline \\*
In the year 2012 Shukla and Mishra we studied Randers conformal change by defining as
\begin{equation}
L(x,y)\rightarrow L^{*}(x, y)=e^{\sigma(x)}L(x, y)+\beta (x, y),
\end{equation}
where $ \sigma(x) $ is a function of x and $\beta(x, y) = b_{i}(x)y^{i}$ is a 1- form on $ M^{n} $. This change generalizes various types of changes. When $ \beta $ = 0, it reduces to a conformal change. When $ \sigma $ = 0, it reduces to a Randers change. When $ \beta=0 $ and $ \sigma $ is a non-zero constant then it reduces to homothetic change.
\newline \\*
In the present paper we have obtained the condition for a Randers conformal change to be projective and find out when a totally geodesic hypersurface $ F^{n-1} $ remains to be a totally geodesic hypersurface $ F^{n-1} $ under the projective Randers conformal change. Further we obatined the condition under which a Finslerian hypersurfaces given by the projective Randers conformal change are projectively flat.
\section{Preliminaries}
\hspace{10pt} Let ($ M^{n} $, L) be a Finsler space $ F^{n} $, where $ M^{n} $ is an n-dimensional differentiable manifold equipped with a fundamental function L. A change in fundamental metric L, defined by equation (1.2), is called Randers conformal change, where $ \sigma(x) $ is conformal factor and function of position only and $\beta(x, y) = b_{i}(x)y^{i}$ is a 1- form on $ M^{n} $. A space equipped with fundamental metric $ L^{*}(x, y) $ is called Randers conformally changed space $ F^{*n} $.\newline\\*
Differentiating equation (1.2) with respect to $ y^{i} $, the normalized supporting element $ l^{*}_{i}=\dot{\partial}_{i}L^{*} $ is given by
\begin{equation}
l^{*}_{i}(x, y) = e^{\sigma(x)}l_{i}(x, y) + b_{i}(x),
\end{equation}
where $ l_{i}=\dot{\partial}_{i}L $ is the normalized supporting element in the Finsler space $ F^{n} $.\newline \\*
Differentiating (2.1) with respect to $ y^{j} $, the angular metric tensor $ h_{ij}^{*}=L^{*}\dot{\partial}_{i}
\dot{\partial}_{j}L^{*} $ is given by
\begin{equation}
h_{ij}^{*}=e^{\sigma(x)}\frac{L^{*}}{L} h_{ij}
\end{equation}
where $ h_{ij}=L\dot{\partial}_{i}\dot{\partial}_{j}L $ is the angular metric tensor in the Finsler space $ F^{n} $ .\newline \\*
Again the fundamental tensor $ g_{ij}^{*}=\dot{\partial}_{i}\dot{\partial}_{j}\frac{L^{*2}}{2}=h_{ij}^{*}+l_{i}^{*}l_{j}^{*} $ is given by
\begin{equation}
g_{ij}^{*}=\tau g_{ij}+b_{i}b_{j}+e^{\sigma(x)}L^{-1}(b_{i}y_{j}+b_{j}y_{i})-\beta e^{\sigma(x)}L^{-3}y_{i}y_{j}
\end{equation}
where we put $y_{i}=g_{ij}(x,y)y^{j}$, $ \tau = e^{\sigma(x)}\frac{L^{*}}{L} $ and $g_{ij}$ is the fundamental tensor of the Finsler space $ F^{n} $. It is easy to see that the det$( g_{ij}^{*})$ does not vanish, and the reciprocal tensor with components $g^{*ij}$ is given by
\begin{equation}
g^{*ij}=\tau^{-1}g^{ij}+\phi y^{i}y^{j}-L^{-1}\tau^{-2}(y^{i}b^{j}+y^{j}b^{i})
\end{equation}
where $ \phi=e^{-2\sigma(x)}(Le^{\sigma(x)}b^{2}+\beta)L^{*-3} $, $ b^{2}=b_{i}b^{i} $, $ b^{i}=g^{ij}b_{j} $ and $ g^{ij} $ is the reciprocal tensor of $ g_{ij} $.\newline \\*
Here it will be more convenient to use the tensors
\begin{equation}
h_{ij}=g_{ij}-L^{-2}y_{i}y_{j}, \;\;\;\;\;\;\;\;\;\;\; a_{i}=\beta L^{-2}y_{i}-b_{i}
\end{equation}
both of which have the following interesting property:
\begin{equation}
h_{ij}y^{j}=0,\;\;\;\;\;\;\;\;\; a_{i}y^{i}=0
\end{equation}
Now differentiating equation (2.3) with respect to $ y^{k} $ and using relation (2.5), the Cartan covariant tensor $C^{*}$ with the components $C_{ijk}^{*}=\dot{\partial}_{k}(\frac{g_{ij}^{*}}{2})$ is given as:
\begin{equation}
C_{ijk}^{*}=\tau [C_{ijk}-\frac{1}{2L^{*}}(h_{ij}a_{k}+h_{jk}a_{i}+h_{ki}a_{j})]
\end{equation}
where $ C_{ijk} $ is (h)hv-torsion tensor of Cartan's connection $ C\Gamma $ of Finsler space $ F^{n} $.\newline \\*
In order to obtain the tensor with the components $ C_{ik}^{*j} $, paying attention to (2.6), we obtain from (2.4) and (2.7),
\begin{eqnarray}
C_{ik}^{*j}=C_{ik}^{j}-\frac{1}{2L^{*}}(h_{i}^{j}a_{k}+h_{k}^{j}a_{i}+h_{ik}a^{j})-(\tau L)^{-1}C_{ikr}y^{j}b^{r}-\\\nonumber \frac{\tau^{-1}}{2LL^{*}}(2a_{i}a_{k}+a^{2}h_{ik})y^{j}
\end{eqnarray}
where $ a_{i}a^{i}=a^{2} $.
We denote by the symbol $ (_{|}) $ the h-covariant differentiation with respect to the Cartan connection $ C\Gamma = (F^{i}_{jk},N^{i}_{j}, C^{i}_{jk}) $ and put
\begin{equation}
2E_{ij}=b_{i|j} + b_{j|i} \;\;\;\;\;\;\; 2F_{ij}=b_{i|j} - b_{j|i}
\end{equation}
Now we deal with well-known functions $ G^{i}(x, y) $ which are (2)p-homogeneous in $ y^{i} $ and are written as $ 2G^{i}=\gamma_{jk}^{i}y^{j}y^{k}  $ by putting $ \gamma_{jk}^{i}=g^{ir}\frac{(\partial_{k}g_{jr}+\partial_{j}g_{kr}-\partial_{r}g_{jk})}{2} $.\newline \\*
Owing to (2.3) and (2.4), a straightforward calculation leads to
\begin{equation}
G^{*i}(x, y)= \frac{\gamma_{jk}^{*i}y^{j}y^{k}}{2}=G^{i}+D^{i}
\end{equation}
where the vector $D^{i}$ is given by\newline
$ D^{i}=\frac{1}{2}\lbrace \tau^{-1}g^{ir}+2\phi y^{i}y^{r}-2L^{(-1)}\tau^{-1}(y^{i}b^{r}+y^{r}b^{i})\rbrace [\tau_{|0}2b_{r}E_{00}+4\beta F_{r0} -2l_{r} -\tau_{|r}+e^{\sigma}L^{(-1)}\lbrace 2\sigma_{|0}(b_{r}+\beta y_{r})-2\beta \sigma_{|r} 4F_{r0}+2E_{00}y_{r}\rbrace -e^{\sigma}L^{(-3)}\lbrace 2\beta_{|0}y_{r}-\beta_{|r}\rbrace - \beta e^{\sigma}L^{(-3)}\lbrace 2\sigma_{|0}y_{r}-\sigma_{|r}\rbrace], $\newline
$F^{i}_{j}=g^{ir}F_{rj}$ and the subscript '0' means the contraction by $y^{i}$.
\section{Relation between projective change and Randers conformal change}
\hspace{10pt} For two Finsler spaces $ F^{n}=(M^{n}, L) $ and $ F^{*n}=(M^{n}, L^{*}) $,if any geodesic on $F^{n}$ is also a geodesic on $F^{*n}$ and the inverse is true, the change $ L\rightarrow L^{*} $ of the
metric is called \textit{projective}. A geodesic on $F^{n}$ is given by a system of differential equations
\begin{equation}
\frac{d^{2}y^{i}}{dt^{2}}+2G^{i}(x, y)=y^{i}, \;\;\;\; y^{i}=\frac{dx^{i}}{dt}
\end{equation}
where $ G^{i}(x, y) $ are (2) p -homogeneous functions in $ y^{i} $. We are now in a position to find a condition for a Randers conformal change to be projective. For this purpose we deal with Euler-Lagrange equations $B_{i}=0$ , where $ B_{i} $ is defined by $ B_{i}=\partial_{i}L-\frac{d(\dot{\partial_{i}}L)}{dt} $. Therefore from the Euler-Lagrange differential equations $ B^{*}_{i}=0 $ for $ F^{*} $ are given by
\begin{center}
$ B^{*}_{i}=e^{\sigma}B_{i}+e^{\sigma}L\partial_{i}\sigma + \partial_{i}\beta - \frac{db_{i}}{dt}=0 $
\end{center}
Thus the above equation can be written as
\begin{equation}
B^{*}_{i}=e^{\sigma}B_{i}+A_{i}
\end{equation}
where $ A_{i} $ is a covariant vector and defined as $ A_{i}=e^{\sigma}L\partial_{i}\sigma + \partial_{i}\beta - \frac{db_{i}}{dt} $.\newline 
Thus we have
\begin{proposition}
Let $F*{n} = (M^{n}, L^{*})$ be an n-dimensional Finsler space obtained from the Randers Conformal change of the Finsler space $F^{n} = (M^{n}, L)$, and 1-form metric, then the Finsler metric $ L^{*} $ is projective if the covariant vector $ A_{i} $ of the equation (3.2) vanishes identically.
\end{proposition}
\section{Hypersurface given by projective Randers conformal change} 
\hspace{10pt} Hereafter, we assume that metrics $ L^{2} $ and $L^{*2} $ are positive-definite respectively and we consider hypersurfaces. According to \cite{M1}, a hypersurface $ M^{n-1} $ of the underlying smooth manifold $ M^{n} $ may be parametrically represented by the equation $ x^{i}=x^{i}(u^{\alpha}), $ where $ u^{\alpha} $ are Gaussian coordinates on $ M^{n-1} $ and Greek indices vary from 1 to n-1. Here we shall assume that the matrix consisting of the projection factors $ B^{i}_{\alpha}=\frac{\partial x^{i}}{\partial u^{\alpha}} $ is of rank n-1. The following notations are also employed:
\begin{center}
$ B^{i}_{\alpha\beta}=\frac{\partial^{2}x^{i}}{\partial u^{\alpha}\partial u^{\beta}}, \;\;\;\;\;\;\;\; B^{i}_{0\beta}=v^{\alpha}B^{i}_{\alpha\beta} $
\end{center}
If the supporting element $ y^{i} $ at a point $ (u^{\alpha}) $ of $ M^{n-1} $ is assumed to be tangential to $ M^{n-1} $, we may then write $ y^{i}=B^{i}_{\alpha}(u)v^{\alpha} $ i.e. $ v^{\alpha} $ is thought of as the supporting element of $ M^{n-1} $ at the point $ (u^{\alpha}) $. Since the function $ \bar{L}(u, v)=L\lbrace x(u), y(u, v)\rbrace $ gives rise to a Finsler metric of $ M^{n-1} $, we get a $ (n-1) $-dimensional Finsler space $ F^{n-1}=\lbrace M^{n-1}, \bar{L}(u, v)\rbrace $. \newline \\*
At each point $ (u^{\alpha}) $ of $ F^{n-1} $, the unit normal vector $ N^{i}(u, v) $ is defined by
\begin{equation}
g_{ij}B^{i}_{\alpha}N^{j}=0, \;\;\;\;\;\; g_{ij}N^{i}N^{j}=1
\end{equation}
If $ B_{i}^{\alpha}, N_{i} $ is the inverse matrix of $ (B^{i}_{\alpha}, N^{i}) $, we have
\begin{center}
$ B^{i}_{\alpha}B_{i}^{\beta}=\delta^{\beta}_{\alpha}, \;\;\;\;\;\;\; B^{i}_{\alpha}N_{i}=0, \;\;\;\;\;\;\; N^{i}N_{i}=1 \;\;\;\;\; $ and $ \;\;\;\;\; B^{i}_{\alpha}B_{j}^{\alpha}+N^{i}N_{j}=\delta^{i}_{j} $.
\end{center}
Making use of the inverse matrix $ (g^{\alpha\beta}) $ of $ (g_{\alpha\beta}) $, we get
\begin{equation}
B^{\alpha}_{i}=g^{\alpha\beta}g_{ij}B^{j}_{\beta},\;\;\;\;\;\;\; N_{i}=g_{ij}N^{j}
\end{equation}
For the induced Cartan's connection $ IC\Gamma = (F^{\alpha}_{\beta\gamma}, N^{\beta}_{\alpha}, C^{\alpha}_{\beta\gamma}) $ on $ F^{n-1} $, the normal curvature vector $ H_{\alpha} $ is given by \cite{M2}
\begin{equation}
H_{\alpha}=N_{i}(B^{i}_{0\beta}+N^{i}_{j}B^{j}_{\beta})
\end{equation}
Consider a Finslerian hypersurface $ F^{n-1}=\lbrace M^{n-1}, \bar{L}(u, v)\rbrace $ of the $ F^{n} $ and another Finslerian hypersurface $ F^{*n-1}=\lbrace M^{n-1}, \bar{L}^{*}(u, v)\rbrace $ of the $ F^{*n} $ given by the Randers conformal change. Let $ N^{i} $ be the unit vector at each point of $ F^{n-1} $ and $ (B^{\alpha}_{i}, N_{i}) $ be the inverse matrix of $ (B^{i}_{\alpha}, N^{i}) $. The function $ B^{i}_{\alpha} $ may be considered as components of $ (n-1) $ linearly independent tangent vectors of $ F^{n-1} $ and they are invariant under Randers conformal change. Thus we shall show that a unit normal vector $ N^{*i}(u, v) $ of $ F^{*n-1} $ is uniquely determined by
\begin{equation}
g^{*}_{ij}B^{i}_{\alpha}N^{*j}=0, \;\;\;\;\;\; g^{*}_{ij}N^{*i}N^{*j}=1
\end{equation} 
Contracting (2.3) by $ N^{i}N^{j} $ and paying attention to (4.1) and the fact that $ l_{i}N^{i}=0 $, we have
\begin{equation}
g^{*}_{ij}N^{i}N^{j}= \tau + (b_{i}N^{i})^{2}
\end{equation}
Therefore we obtain
\begin{center}
$ g^{*}_{ij}\lbrace \pm \frac{N^{i}}{\sqrt{\tau + (b_{i}N^{i})^{2}}}\rbrace \lbrace \pm \frac{N^{j}}{\sqrt{\tau + (b_{i}N^{i})^{2}}}\rbrace =1 $
\end{center}
Hence we can put
\begin{equation}
N^{*i}= \frac{N^{i}}{\sqrt{\tau + (b_{i}N^{i})^{2}}}
\end{equation}
where we have chosen the positive sign in order to fix an orientation. \newline \\*
Using equation (2.3), (4.6) and from first condition of (4.4) we have
\begin{equation}
(b_{i}B^{i}_{\alpha}+e^{\sigma(x)}l_{i}B^{i}_{\alpha})\frac{b_{j}N^{j}}{\sqrt{\tau + (b_{i}N^{i})^{2}}}=0
\end{equation}
If $ b_{i}B^{i}_{\alpha}+e^{\sigma(x)}l_{i}B^{i}_{\alpha}=0 $, then contracting it by $ v^{\alpha} $ and using $ y^{i}=B^{i}_{\alpha}v^{\alpha} $ we get $ \beta + e^{\sigma(x)}L=L^{*}=0 $ which is contradiction to the assumption that $ L^{*}>0 $. Hence $ b_{i}N^{i}=0 $. Therefore equation (4.6) can be written as
\begin{equation}
N^{*i}=\frac{1}{\sqrt{\tau}}N^{i}
\end{equation}
Summarizing the above Shukla, Chaubey and Mishra \cite{SCM} obtained the following result
\begin{proposition}
If $ \lbrace (B^{i}_{\alpha}, N^{i}), \; \alpha=1,2,...(n-1) \rbrace $ be the filed of linear frame of the Finsler space $ F^{n} $, there exist a field of linear frame $ \lbrace (B^{i}_{\alpha}, N^{*i}=\frac{1}{\sqrt{\tau}}N^{i}), \; \alpha=1,2,...(n-1)\rbrace $ of the Finsler space $ F^{*n} $ such that (15) is satisfied along $ F^{*n-1} $ and then $ b_{i} $ is tangential to both the hypersurfaces $ F^{n-1} $ and $ F^{*n-1} $.
\end{proposition}
The quantities $ B_{i}^{*\alpha} $ are uniquely defined along $ F^{*n-1} $ by
\begin{center}
$  B_{i}^{*\alpha}=g^{*\alpha\beta}g^{*}_{ij}B^{j}_{\beta} $  
\end{center}
where $ g^{*\alpha\beta} $ is the inverse matrix of $ g^{*}_{\alpha\beta} $. Let $ (B_{i}^{*\alpha}, N^{*}_{i}) $ be the inverse matrix of $ (B^{i}_{\alpha}, N^{*i}) $, then we have  
\begin{center}
$ B^{i}_{\alpha}B_{i}^{*\beta}=\delta^{\beta}_{\alpha},\;\;\;\;\;  B^{i}_{\alpha}N^{*}_{i}=0, \;\;\;\;\;\;N^{*i}N^{*}_{i}=1 $
\end{center} 
Furthermore $ B^{i}_{\alpha}B_{j}^{*\alpha}+ N^{*i}N^{*}_{j}=\delta^{i}_{j}$. We also get $ N^{*}_{i}=g^{*}_{ij}N^{*j} $ which in view of (2.1), (2.3) and (4.8) gives
\begin{equation}
N^{*}_{i}=\sqrt{\tau}N_{i}
\end{equation}
Now we assume that a Randers conformal change of the metric is projective. \newline
Using (2.10) and Proposition 3.1, we have
\begin{equation}
D^{i}=G^{*i}-G^{i}
\end{equation}
Since $ D^{i}_{j}=\dot{\partial}_{j}D^{i} $ and $ N^{i}_{j}=\dot{\partial}_{j}G^{i} $ the above gives
\begin{equation}
D^{i}_{j}=N^{*i}_{j}-N^{i}_{j}
\end{equation}
Further contracting the above equation by $ N_{i}B_{\alpha}^{j} $ we have
\begin{equation}
N_{i}D^{i}_{j}B^{j}=0
\end{equation}
If each geodesic of $F^{(n−1)}$ with respect to the induced metric is also a geodesic of $F^{(n)}$ , then $F^{(n−1)}$ is called totally geodesic. A totally geodesic hypersurface $F^{(n−1)}$ is characterized by $ H_{\alpha}=0 $. \newline
From equation (4.3), (4.9) and (4.11) we have
\begin{equation}
H^{*}_{\alpha}=\sqrt{\tau}H_{\alpha}+N_{i}D^{i}_{j}B^{j}
\end{equation}
Thus using the equation (4.12) in the above equation we have
\begin{equation}
H^{*}_{\alpha}=\sqrt{\tau}H_{\alpha}
\end{equation} 
Thus we have
\begin{theorem}
A hypersurface $F^{n - 1}$ of a Finsler space $ F^{n} (n>3) $ is totally geodesic, if and only if the hypersurface $F^{*(n - 1)}$ of the space $F^{*n}$ obtained from $F^{n}$ by a projective Randers conformal change, is totally geodesic.
\end{theorem}
\section{Hypersurfaces of Projectively Flat Finsler spaces}
In this section, we shall consider a projective Randers conformal change and we are concerned with the Berwald connection $ B\Gamma $ on $ F^{n}=(M^{n}, L) $ and $ B\bar{\Gamma} $ on $ F^{*n}=(M^{*n}, L) $. In the theory of projective changes in Finsler spaces, we have two essential projective invariants, one is the Weyl torsion tensor $ W_{ij}^{h} $ and the other is the Douglas tensor $ D^{h}_{ijk} $, so that under the projective Randers conformal change, we get $ W_{ij}^{*h} $=$ W_{ij}^{h} $ and $ D^{*h}_{ijk} $=$ D^{h}_{ijk} $. \newline \\*
Now we are concerned with a projectively flat Finsler space defined as
follows: If there exists a projective change $ L\rightarrow L^{*} $ of a Finsler space $ F^{n}=(M^{n}, L) $ such that the Finsler space $ F^{*n}=(M^{*n}, L) $ is a locally Minkowski space then $ F^{n} $ is called projectively flat Finsler space. We have already known the following theorems: 
\begin{theorem}
\cite{K2} A Finsler space $ F^{n} (n>2) $ is projectively flat, if and only if $ W_{ij}^{h}=0 $ and $ D^{h}_{ijk}=0 $.
\end{theorem}
\begin{theorem}
\cite{SSA} A Finsler space $ F^{n} (n>3) $ is projectively flat then the
totally geodesic hypersurface $F^{n - 1}$ is also projectively flat.
\end{theorem}
Thus from theorem (4.1), theorem (5.1) and theorem (5.2) we have
\begin{theorem}
Let $ F^{n} (n>3) $ be a projectively flat Finsler space. If the hypersurface $F^{n - 1}$ is totally geodesic, then the hypersurface $F^{*(n - 1)}$ of the space $F^{*n}$ obtained from $F^{n}$ by a projective Randers conformal change, is projectively flat.
\end{theorem}

\end{document}